\documentclass[12pt]{article}
\usepackage[T2A]{fontenc}
\usepackage[utf8]{inputenc}
\usepackage{amssymb,amsmath,amsfonts,amsthm,amscd,latexsym,indentfirst}
\usepackage{epsfig,geometry}

\def\id{\textrm{id}}

\def\rank{\textrm{rank}\,}

\def\Span{\mathrm{Span}}
\def\Aut{\mathrm{Aut}}

\begin{document}

\begin{center}
{\Large
Rota---Baxter operators on a sum of fields}

V. Gubarev
\end{center}

\begin{abstract}
We count the number of all Rota---Baxter operators on a finite direct sum 
$A = F\oplus F\oplus \ldots \oplus F$ of fields
and count all of them up to conjugation with an automorphism.
We also study Rota---Baxter operators on $A$ corresponding to a decomposition of $A$ 
into a~direct vector space sum of two subalgebras.
We show that every algebra structure induced on $A$ by a Rota---Baxter of nonzero weight is isomorphic to $A$.

\end{abstract}

\section{Introduction}

Given an algebra $A$ and a scalar $\lambda\in F$, where $F$ is a~ground field, 
a~linear operator $R\colon A\rightarrow A$ is called a Rota---Baxter operator 
(RB-operator, for short) on $A$ of weight~$\lambda$ if the following identity
\begin{equation}\label{RB}
R(x)R(y) = R( R(x)y + xR(y) + \lambda xy )
\end{equation}
holds for any $x,y\in A$. The algebra $A$ is called Rota---Baxter algebra (RB-algebra).

G. Baxter in 1960 introduced the notion of Rota---Baxter operator~\cite{Baxter} 
as natural generalization of by parts integration formula. In~1960--1970s
such operators were studied by G.-C.Rota~\cite{Rota}, P.~Cartier~\cite{Cartier},
J. Miller~\cite{Miller}, F. Atkinson~\cite{Atkinson} and others.

In 1980s, the deep connection between constant solutions 
of the classical Yang---Baxter equation from mathematical 
physics and RB-operators on a semi\-simple finite-dimensional Lie algebra 
was discovered by A. Belavin and V. Drinfel'd~\cite{BelaDrin82}
and M. Semenov-Tyan-Shanskii~\cite{Semenov83}. 

About different connections of Rota---Baxter operators 
with symmetric polynomials, quantum field renormalization, Loday algebras,
shuffle algebra see in the monograph~\cite{GuoMonograph} 
written by L. Guo in 2012.

In the paper, we study Rota---Baxter operators on a finite direct sum 
$A = F\oplus F\oplus \ldots \oplus F$ of $n$ copies of a field~$F$.
We continue investigations fulfilled by S. de Bragan\c{c}a in 1975~\cite{Braga}
and by H. An and C. Bai in 2008~\cite{AnBai}.
Since all RB-operators on $A$ of weight zero are trivial~\cite{Unital}, i.e., equal to~0,
we study only RB-operators on $A$ of nonzero weight~$\lambda$.

In~\S2, we formulate some preliminaries about RB-operators,
including splitting RB-operators which are projections on a subalgebra~$A_1$ 
parallel to another one $A_2$ provided the direct vector space sum decomposition $A = A_1\dot+A_2$.

In~\S3, we show that RB-operators on $A$ of nonzero weight~$\lambda$
are in bijection with 2-colored transitive subtree acyclic digraphs
(subtree acyclic digraphs were defined by F.~Harary et al. in 1992~\cite{Harary})
or equivalently with labeled rooted trees on $n+1$ vertices
with 2-colored non-root vertices. For the last, we apply 
the result of R. Castelo and A.~Siebes~\cite{Castelo}.
Thus, the number of all RB-operators on $A$ of nonzero weight~$\lambda$ 
equals $2^n(n+1)^{n-1}$. 
With the help of the bijection, we show that splitting RB-operators on~$A$ 
of nonzero weight~$\lambda$ are in one-to-one correspondence
with labeled rooted trees on $n+1$ vertices with properly 2-colored non-root vertices.
We also study the number of all RB-operators and all splitting RB-operators on~$A$
up to conjugation with an automorphism of~$A$.

In 2012, D. Burde et al. initiated to study so called post-Lie algebra structures~\cite{BurdeDekimpeVercammen12}.
One of the questions arisen in the area~\cite{BurdeDekimpeVercammen12,BurdeDekimpe13,BurdeGub18} is the following one: 
starting with a~semisimple Lie algebra endowed RB-operator of weight~1 
what kind of Lie algebras we will get under the new Lie bracket $[R(x),y] + [x,R(y)] + [x,y]$?
Such problems could be stated not only for Lie algebras but also for associative or commutative ones.
In~\S4, we show that every algebra structure induced on 
a finite direct sum $A$ of fields by a~Rota---Baxter operator 
of nonzero weight is isomorphic to $A$ itself.

\section{Preliminaries}

Trivial RB-operators of weight $\lambda$ are zero operator and $-\lambda\id$.

{\bf Statement 1}~\cite{GuoMonograph}.
Given an RB-operator $R$ of weight $\lambda$,

a) the operator $-R-\lambda\id$ is an RB-operator of weight $\lambda$,

b) the operator $\lambda^{-1}R$ is an RB-operator of weight 1, provided $\lambda\neq0$.

Given an algebra $A$, let us define a map $\phi$ on the set of all RB-operators on $A$
as $\phi(R)=-R-\lambda(R)\id$. It is clear that $\phi^2$ coincides with the identity map.

{\bf Statement 2}~\cite{BGP}.
Given an algebra $A$, an RB-operator $R$ on $A$ of weight $\lambda$,
and $\psi\in\Aut(A)$, the operator $R^{(\psi)} = \psi^{-1}R\psi$
is an RB-operator on $A$ of weight $\lambda$.

{\bf Statement 3}~\cite{GuoMonograph}.
Let an algebra $A$ to split as a vector space
into the direct sum of two subalgebras $A_1$ and $A_2$.
An operator $R$ defined as
\begin{equation}\label{Split}
R(a_1 + a_2) = -\lambda a_2,\quad a_1\in A_1,\ a_2\in A_2,
\end{equation}
is RB-operator on $A$ of weight $\lambda$.

Let us call an RB-operator from Statement~3 as
{\it splitting} RB-operator with subalgebras $A_1,A_2$.
Note that the set of all splitting RB-operators on 
an algebra $A$ is in bijection with all decompositions 
$A$ into a direct sum of two subalgebras $A_1,A_2$.

{\bf Remark~1}.
Given an algebra $A$, let $R$ be a splitting RB-operator
on $A$ of weight $\lambda$ with subalgebras $A_1,A_2$.
Hence, $\phi(R)$ is an RB-operator of weight $\lambda$ and
$$
\phi(R)(a_1+a_2) = -\lambda a_1,\quad a_1\in A_1,\ a_2\in A_2.
$$
So $\phi(R)$ is splitting RB-operator with the same subalgebras $A_1, A_2$.

{\bf Lemma~1}~\cite{BGP}.
Let $A$ be a unital algebra, $R$ be an RB-operator on $A$ of nonzero weight~$\lambda$.
If $R(1)\in F$, then $R$ is splitting.

We call an RB-operator $R$ satisfying the conditions of Lemma~1 as 
{\it inner-splitting} one.

\newpage

{\bf Lemma~2}~\cite{Unital}.
Let $A = A_1\oplus A_2$ be an algebra, 
$R$ be an RB-operator on $A$ of weight~$\lambda$.
Then the induced linear map 
$P\colon A_1 \to A_1$ defined by the formula 
$P(x_1+x_2) = \mathrm{Pr}_{A_1}(R(x_1))$, $x_1\in A_1$, $x_2\in A_2$,
is an RB-operator on $A_1$ of weight~$\lambda$.

\section{RB-operators on a sum of fields}

{\bf Statement~4}~\cite{AnBai,Braga,Unital}.
Let $A = Fe_1\oplus Fe_2\oplus\ldots\oplus Fe_n$
be a direct sum of copies of a~field~$F$.
A linear operator $R(e_i) = \sum\limits_{k=1}^n r_{ik}e_k$,
$r_{ik}\in F$, is an RB-operator on $A$ of weight~1
if and only if the following conditions are satisfied:

(SF1) $r_{ii} = 0$ and $r_{ik}\in\{0,1\}$
or $r_{ii} = -1$ and $r_{ik}\in\{0,-1\}$ for all $k\neq i$;

(SF2) if $r_{ik} = r_{ki} = 0$ for $i\neq k$,
then $r_{il}r_{kl} = 0$ for all $l\not\in\{i,k\}$;

(SF3) if $r_{ik}\neq0$ for $i\neq k$,
then $r_{ki} = 0$ and
$r_{kl} = 0$ or $r_{il} = r_{ik}$ for all $l\not\in\{i,k\}$.

{\bf Example}~\cite{Atkinson,Miller}.
The following operator is an RB-operator on $A$ of weight~1:
$$
R(e_i) = \sum\limits_{l=i+1}^s e_l,\ 1\leq i<s,\quad
R(e_s) = 0,\quad  
R(e_i) = -\sum\limits_{l=i}^{n} e_l,\ s+1\leq i\leq n.
$$

{\bf Remark~2}.
It follows from (SF3) that $r_{ik}r_{ki} = 0$ for all $i\neq k$.
In \cite{AnBai}, the statement of Statement~4 was formulated with this equality and (SF1)
but without (SF2) and the general version of (SF3).
That's why the formulation in \cite{AnBai} seems to be not complete.

{\bf Remark~3}.
The sum of fields in Statement~4 can be infinite.

In advance, we will identify an RB-operator on~$A$ with its matrix.

Let us calculate the number of different RB-operators of nonzero weight $\lambda$ on
$A = Fe_1\oplus Fe_2\oplus\ldots\oplus Fe_n$. By Statement 1a, we may assume that $\lambda = 1$.
For $n = 1$, we have only two RB-operators $\{0,-\id\}$. 
For $n = 2$ we have 12 cases~\cite{AnBai}:
\begin{gather*}
\begin{pmatrix}
0 & 0 \\
0 & 0 \\
\end{pmatrix},\
\begin{pmatrix}
-1 & 0 \\
0 & -1 \\
\end{pmatrix},\
\begin{pmatrix}
0 & 0 \\
1 & 0 \\
\end{pmatrix},\
\begin{pmatrix}
-1 & 0 \\
-1 & -1 \\
\end{pmatrix},\
\begin{pmatrix}
-1 & 0 \\
0 & 0 \\
\end{pmatrix},\
\begin{pmatrix}
0 & 0 \\
0 & -1 \\
\end{pmatrix},\\
\begin{pmatrix}
0 & 0 \\
-1 & -1 \\
\end{pmatrix},\
\begin{pmatrix}
-1 & 0 \\
1 & 0 \\
\end{pmatrix},\
\begin{pmatrix}
-1 & -1 \\
0 & 0 \\
\end{pmatrix},\
\begin{pmatrix}
0 & 1 \\
0 & -1 \\
\end{pmatrix},\
\begin{pmatrix}
0 & 1 \\
0 & 0 \\
\end{pmatrix},\
\begin{pmatrix}
-1 & -1 \\
0 & -1 \\
\end{pmatrix}.
\end{gather*}
Here we identify an RB-operator with its matrix $R\in M_2(F)$
by the rule $R(e_i) = \sum\limits_{k=1}^n r_{ik}e_k$.

For $n = 3$, we have $8\cdot16 = 128$ variants~\cite{AnBai}:
\begin{gather*}
\begin{pmatrix}
a&0&0\\
0&b&0\\
0&0&c
\end{pmatrix},\
\begin{pmatrix}
a&0&0\\
0&b&0\\
2c+1&2c+1&c
\end{pmatrix},\
\begin{pmatrix}
a&0&0\\
0&b&0\\
2c+1&0&c
\end{pmatrix},\
\begin{pmatrix}
a&0&0\\
0&b&0\\
0&2c+1&c
\end{pmatrix},\\
\allowdisplaybreaks
\begin{pmatrix}
a&0&0\\
2b+1&b&0\\
0&0&c
\end{pmatrix},\
\begin{pmatrix}
a&0&0\\
2b+1&b&0\\
2c+1&2c+1&c
\end{pmatrix},\
\begin{pmatrix}
a&2a+1&0\\
0&b&0\\
0&0&c
\end{pmatrix},\
\begin{pmatrix}
a&2a+1&0\\
0&b&0\\
2c+1&2c+1&c
\end{pmatrix},\\
\begin{pmatrix}
a&2a+1&2a+1\\
0&b&0\\
0&0&c
\end{pmatrix}\!,\
\begin{pmatrix}
a&2a+1&2a+1\\
0&b&0\\
0&2c+1&c
\end{pmatrix}\!,\
\begin{pmatrix}
a&0&0\\
2b+1&b&2b+1\\
0&0&c
\end{pmatrix}\!,\
\begin{pmatrix}
a&0&0\\
0&b&2b+1\\
0&0&c
\end{pmatrix}\!,\\
\begin{pmatrix}
a&0&0\\
2b+1&b&2b+1\\
2c+1&0&c
\end{pmatrix},\
\begin{pmatrix}
a&2a+1&2a+1\\
0&b&2b+1\\
0&0&c
\end{pmatrix},\
\begin{pmatrix}
a&0&2a+1\\
0&b&0\\
0&0&c
\end{pmatrix},\
\begin{pmatrix}
a&0&2a+1\\
2b+1&b&2b+1\\
0&0&c
\end{pmatrix}
\end{gather*}
for $a=r_{11},b=r_{22},c=r_{33}\in\{0,-1\}$.

For $n = 4$, computer can help to state that there are exactly 2000 RB-operators
of weight~1 on $A$. Thus, we get the first four terms from the sequence A097629~\cite{OEIS}.

{\bf Theorem 1}. 
Let $A = Fe_1\oplus Fe_2\oplus\ldots\oplus Fe_n$
be a direct sum of copies of a~field~$F$.
The number of different RB-operators on $A$ of nonzero weight~$\lambda$
equals $2^n(n+1)^{n-1}$.

{\sc Proof}.
Let $R$ be an RB-operator on $A$ of weight~$\lambda$.
We may assume that $\lambda = 1$. We follow the previous notations. 
We have $2^n$ variants to choose the values of the elements $r_{ii}$, $i=1,\ldots,n$.
The choice of any of them, say $r_{ii}$, influences only on the possible signs of all elements $r_{ik}$, $k\neq i$.
So, we may put $r_{ii} = 0$ for all $i$ and fix the factor $2^n$ for the answer.

Now, we want to construct a directed graph $G$ on $n$ vertices by any matrix $R = (r_{ij})_{i,j=1}^n$ with chosen $r_{ii} = 0$. 
We consider the matrix $R$ as the adjacency matrix of a~directed graph $G$. Let us interpretate conditions (SF2) and (SF3) in terms of digraphs.
Firstly, we rewrite (SF3) as two conditions:

(SF3a) if $r_{ik}\neq0$ for $i\neq k$, then $r_{ki} = 0$;

(SF3b) if $r_{ik}\neq0$ for $i\neq k$, then $r_{kl} = 0$ or $r_{il} = r_{ik}$ for all $l\not\in\{i,k\}$.

The condition (SF3a) says that if we have an edge between two vertices $i\neq k$, then the direction of such edge is well-defined,
so, it is a correctness of getting a digraph by the matrix $R$. 
In graph theory, the condition (SF3b) is called {\it transitivity}, i.e., if have edges $(i,k)\in E$ and $(k,l)\in E$,
then we have an edge $(i,l)\in E$.

Secondly, we read the condition (SF2) in terms of digraphs in such way: 
there are no in~$G$ induced subgraphs isomorphic to $H$ with $V(H) = \{i,k,l\}$ and $E(H) = \{(i,l),(k,l)\}$ (see Pict.~1). 
In~\cite{Castelo} the subgraph~$H$ was called {\it immorality}, thus,
a digraph without immoralities is called {\it moral} digraph~\cite{Lauritzen}.

\begin{figure}[ht]
\centering
\includegraphics[height = 2.3cm]{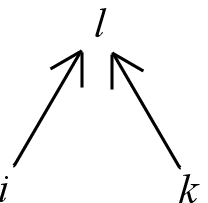} \\
\smallskip
{\sc Picture~1}. The forbidden induced subgraph $H$ on three vertices $\{i,k,l\}$ due to (SF2)
\end{figure}

We may reformulate our problem of counting the number $N$ of different RB-operators on $A$ of nonzero weight~$\lambda$
in such way: What is the number of all transitive moral transitive digraphs on $n$ vertices? 
In terms of~\cite{Castelo}, the last is the same as the number of all moral TDAGs on $n$ vertices, 
here TDAG is the abbreviation for Transitive Directed Acyclic Graph (we are interested on 
transitive digraphs which are surely acyclic).
In the graph-theoretic context, moral DAGs are known as subtree acyclic digraphs~\cite{Harary}.
Thus, 
\begin{multline}
N/2^n = \# \{\mbox{moral TDAGs on } n\mbox{ vertices}\} \\
  = \# \{\mbox{transitive subtree acyclic digraphs on } n\mbox{ vertices}\}.
\end{multline}  

In~\cite{Castelo}, the authors constructed a bijection between the set of moral TDAGs 
on $n$~vertices and the set of labeled rooted trees on $n+1$ vertices as follows (see Pict.~2). 
Define the function $f(i)$ for a vertex $i$ by induction. 
For a source $i$ (i.e., such a vertex~$i$ that there are no edges $(j,i)$ in a digraph), we put $f(i) = 0$.
For a not-source vertex $j$, we may find the unique source $i$ such that there exists 
a directed path $p$ from $i$ to $j$. So, we define $f(j)$ as the length of~$p$. 
Now, we construct a labeled rooted tree $T = (U,F)$ by a~moral TDAG $G = G(V,E)$:
$$
U = V\cup\{0\},\quad 
F = \{(0,i)\mid f(i) = 0\}\cup \{(i,j)\mid (i,j)\in E,\,f(i)=f(j)-1\}.
$$

\begin{figure}[ht]
\centering
\includegraphics[height = 5cm]{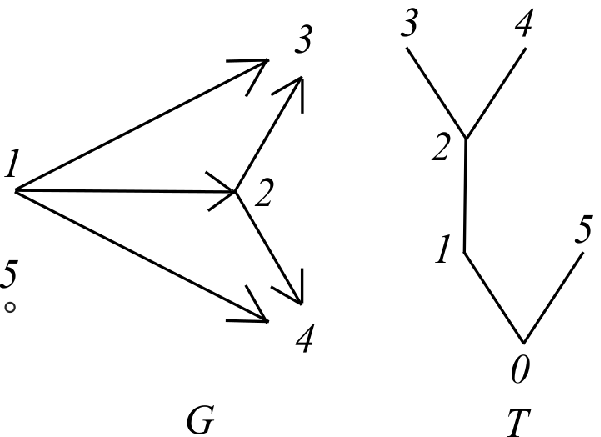} \\
\smallskip
{\sc Picture~2}. The corresponding graph $G$ and tree $T$ to the RB-operator 
$R(e_1) = e_2 + e_3 + e_4$, 
$R(e_2) = -e_2 - e_3 - e_4$,
$R(e_3) = -e_3$, 
$R(e_4) = 0$, 
$R(e_5) = -e_5$.
\end{figure}

Applying the above constructed correspondence, the number of moral TDAGs on $n$~vertices
equals $(n+1)^{n-1}$ by the Cayley theorem, and so $N = 2^n(n+1)^{n-1}$.
Theorem is proved.

Below we will apply the easy fact that $\mathrm{Aut}(A)\cong S_n$.
It could be derived, e.g., from the Molin---Wedderburn---Artin theory, in particular 
from the uniqueness up to a~rearrangement of summands
of decomposition of a semisimple finite-dimensional associative algebra into a~finite direct sum of simple ones.

{\bf Corollary 1}~\cite{Braga}. 
Let $A = Fe_1\oplus Fe_2\oplus\ldots\oplus Fe_n$
be a direct sum of copies of a~field~$F$ and $R$ be an RB-operator on $A$ of nonzero weight~1.
There exists an automorphism~$\psi$ of~$A$ such that the matrix of the operator $R^{(\psi)}$ 
in the basis $e_1,\ldots,e_n$ is an upper-triangular matrix with entries $r_{ij}\in\{0,\pm1\}$ and $r_{ii}\in \{0,-1\}$. 

{\sc Proof}.
As we did in the proof of Theorem~1, we define by~$R$ a labeled rooted tree~$T$.
Define $t = \max\{f(i)\mid i\in V(T)\}$ and $k_j = \# \{i\mid f(i) = j\}$.
We may reorder indexes $1,2,\ldots,n$ by action of a permutation from $S_n\cong \mathrm{Aut}(A)$ in a way such that 
\begin{gather*}
f(1) = \ldots = f(k_0) = 0,\\
f(k_0+1) = \ldots = f(k_0 + k_1) = 1,\\ 
\ldots \\ 
f(n-k_t+1) = \ldots = f(n) = t.
\end{gather*}
Due to the definition of $T$, we get the upper-triangular matrix. 
The restrictions on the values of elements immediately 
follow from Statement~4.

{\bf Corollary 2}. 
There is a bijection between the set of RB-operators of nonzero weight~$\lambda$ on 
$Fe_1\oplus Fe_2\oplus\ldots\oplus Fe_n$ and 

a) the set of 2-colored subtree acyclic digraphs on $n$ vertices;

b) the set of labeled rooted trees on $n+1$ vertices with 2-colored non-root vertices.

Now, we want to compute the number $r_n$ of RB-operators of nonzero weight $\lambda$ 
on $A = Fe_1\oplus\ldots \oplus Fe_n$ which lie in different orbits under the action 
of the automorphism group $\mathrm{Aut}(A) \cong S_n$. 
The group $\mathrm{Aut}(A)$ acts on the set of RB-operators of weight~$\lambda$
in the way described in Statement~2,
$\psi\colon R\to R^{(\psi)} = \psi^{-1}R\psi$.

In a light of Corollary~2b, we may interpretate the number $r_n$ as the number
of {\bf unlabeled} rooted trees on $n+1$ vertices with 2-colored non-root vertices.
It is exactly the sequence A000151~\cite{OEIS}, the first eight values are
2, 7, 26, 107, 458, 2058, 9498, 44947 etc.
Let us fix that in advance we will use two colors: white and black,
white color corresponds to the case $r_{ii} = 0$ and black color corresponds to $r_{ii} = -\lambda$.
Considering the rooted tree $T$ with $n+1$ vertices, we may assume that
the root is colored in the third color, say grey.

Note that the map $\phi$ acts on a labeled (or unlabeled) rooted tree $T$ on $n+1$ vertices
with 2-colored non-root vertices as follows. The $\phi$ interchanges a color in every non-root vertex.

Let us describe splitting RB-operators of nonzero weight $\lambda$ on $A$.

{\bf Theorem 2}. 
An RB-operator $R$ of nonzero weight $\lambda$ on $A = Fe_1\oplus\ldots\oplus Fe_n$ 
is splitting if and only if the corresponding (labeled) rooted tree $T = T(R)$ 
on $n+1$ vertices is properly colored. 

{\sc Proof}.
Wuthout loss of generality, we put $\lambda = 1$.
For simplicity, let us consider the graph $T' = T\setminus\{\mbox{root}\}$,
which is a forest in general case. 

Let us prove the statement by induction on $n$. 
For $n = 1$, we have either $R = 0$ 
(the only non-root vertex is white) or $R = -\lambda\id$ (the only non-root vertex is black), 
both RB-operators are splitting with subalgebras $F$ and $(0)$.

Suppose that we have proved Theorem~2 for all natural numbers less than $n$.
Let a~graph $T'$ with $n$ vertices be disconnected, denote by $T_1,\ldots,T_k$ the connected components of $T'$.
So, $A = A_1\oplus \ldots A_k$ for $A_s = \Span\{e_j\mid j\in V(T_s)\}$. 
Define $R_s$ as the induced RB-operator $R|_{A_s}$ (see Lemma~2). 
By the definition, $R$ is splitting if and only if $A = \ker(R)\dot+\ker(R+\id)$
or equivalently $A_s = \ker(R_s)\dot+\ker(R_s+\id)$, $s=1,\ldots,k$. 
By the induction hypothesis, we have such decomposition for every $s$
if and only if the coloring of $T_s$ is proper.

Now consider the case when $T'$ is connected. 
We may assume that $e_1$ corresponds to the vertex~1, the only source in $G$,
and $\{2,\ldots,k\}$ is the set of all vertices of $G$ 
with the value of $f(x)$ equal to~1.
We also define $T_s$ for $s=2,\ldots,k$ as the connected component 
of $T'\setminus\{1\}$ which contains the vertex~$s$.
Note that $R$ induces the RB-operator of weight~$\lambda$
on the subalgebra $A_s = \Span\{e_j\mid j\in V(T_s)\}$ for all $s$ by Lemma~2.

The condition of $R$ to be splitting is equivalent to the condition
\begin{equation}\label{rankCond}
\rank(R) + \rank(R+\id) = n. 
\end{equation}
Analysing the $e_1$-coordinate, we have
$$
n = \rank(R) + \rank(R+\id) \geq 1 + \rank(R') + \rank(R'+\id) 
$$
for $R'$, the induced RB-operator on the subalgebra $\Span\{e_j\mid j\geq2\}$.
Thus, $\rank(R') + \rank(R'+\id) = n-1$, i.e. $R'$ is splitting or equivalently
$R|_{A_s}$ is spplitting for every $s=2,\ldots,k$. 
By the induction hypothesis, the graph $T'\setminus\{1\}$ is properly 2-colored.
It remains to prove that the vertices $2,\ldots,k$ are colored in the same color
and the vertex~1 is colored in another one.

Up to the action of $\phi$, which preserves the splitting structure of an RB-operator (see Remark~1),
we may assume that the vertex~1 is colored in white. 
Since we know that $\rank(R+\id) = \rank(R'+\id) + 1$,
we have to state the equality $\rank(R) = \rank(R')$. 
So, the condition~\eqref{rankCond} is fulfilled if and only if the first 
row $(0,1,1,\ldots,1)$ of the matrix~$R$ is linearly expressed via other rows. 
By the definition of the matrix $R$, the vertices $2,\ldots,k$ have to be colored in black.
Theorem is proved.

{\bf Corollary 3}.
An RB-operator $R$ of nonzero weight $\lambda$ on $A = Fe_1\oplus\ldots\oplus Fe_n$ 
is inner-splitting if and only if in $T = T(R)$ all vertices with even value of $f$
are colored in one color and all vertices with odd value of $f$ are colored in another color.

{\sc Proof}.
Up to $\phi$, we may assume that $R(1) = 0$. 
Thus, any vertex with the value of $f(x)$ equal to~0 has to be colored in white.
By Theorem~2, $T' = T\setminus\{\mbox{root}\}$ is properly 2-colored,
so, all vertices with the value of $f(x)$ equal to~1 are colored in black,
all vertices with the value of $f(x)$ equal to~2 are colored in white and so on.

Now, we collect all our knowledges about all RB-operators (in Table~1) 
and all nonisomorphic RB-operators (in Table~2) 
of nonzero weight on a sum of fields $A = Fe_1\oplus Fe_2\oplus \ldots\oplus Fe_n$. 

We have noticed that the first values of number of splitting RB-operators 
coincides with the sequence A007830~\cite{OEIS} (in labeled case) 
and coincides with the sequence A000106~\cite{OEIS} (in unlabeled case). 
Actually it should be proven for all $n$.

{\bf Remark~4}. 
Counting rooted trees on $n+1$ vertices with properly 2-colored non-root vertices
is not the same as counting properly 2-colored forests on $n$ vertices. 

\newpage

\begin{center}
{\sc Table 1}. 
Number of RB-operators of nonzero weight on a sum of $n$ fields

\medskip

\begin{tabular}{c|c|c|l}
Class of & Description & formula and & \!first \\
RB-operators & & OEIS~\cite{OEIS} & \!5 values  \\
\hline
all & labeled rooted trees on $n+1$ vertices & $2^n(n+1)^{n-1}$  & \!$2, 12, 128$, \\ 
    & with 2-colored non-root vertices & A097629 & \!$2000, 41472$\! \\
\hline
splitting & labeled rooted trees on $n+1$ vertices & \!$2(n+2)^{n-1}$\,?!\! & \!$2, 8, 50$, \\ 
    & with properly 2-colored non-root vertices & A007830\,?! & \!$432,4802$ \\
\hline
\!inner-splitting\! & labeled rooted trees on $n+1$ vertices & $2(n+1)^{n-1}$ & \!$2, 6, 32$, \\ 
                                                               & (twice) & $2\cdot$A000272 & \!$250,2592$ \\
\hline
non-splitting & \!labeled rooted trees on $n+1$ vertices with\! & --- & \!$0, 4, 78$, \\ 
    & improperly 2-colored non-root vertices & & \!$1568, 36670$\! \\
\end{tabular}
\end{center}

\medskip

\begin{center}
{\sc Table 2}. 
Number of RB-operators of nonzero weight on a sum of $n$ fields \\
(up to conjugation with an automorphism)

\medskip

\begin{tabular}{c|c|c|l}
Class of & Description & OEIS~\cite{OEIS} & first 5 values \\
RB-operators & & & \\
\hline
all & rooted trees on $n+1$ vertices & A000151 & $2, 7, 26, 107, 458$ \\
    & with 2-colored non-root vertices & & \\
\hline
splitting & rooted trees on $n+1$ vertices with & A000106\,?! & $2, 5, 12, 30, 74$ \\
    & properly 2-colored non-root vertices & & \\
\hline
inner-splitting & rooted trees on $n+1$ vertices (twice) & $2\cdot$A000081 & $2, 4, 8, 18, 40$ \\
\hline
non-splitting & rooted trees on $n+1$ vertices with & --- & $0, 2, 14, 77, 384$ \\
    & improperly 2-colored non-root vertices & & \\
\end{tabular}
\end{center}

\medskip

Let us write down all non-splitting pairwise nonisomorphic RB-operators for $n = 2,3$.

{\bf Statement 5}.
Up to $\phi$, we have the following non-splitting pairwise nonisomorphic RB-operators 

a) for $n=2$: $R(e_1) = e_2$, $R(e_2) = 0$;

b) for $n=3$:

(RB1) $R(e_1) = e_2 + e_3$, $R(e_2) = e_3$, $R(e_3) = 0$,

(RB2) $R(e_1) = e_2 + e_3$, $R(e_2) = e_3$, $R(e_3) = -e_3$,

(RB3) $R(e_1) = e_2 + e_3$, $R(e_2) = -e_2-e_3$, $R(e_3) = -e_3$,

(RB4) $R(e_1) = e_2 + e_3$, $R(e_2) = R(e_3) = 0$,

(RB5) $R(e_1) = e_2 + e_3$, $R(e_2) = -e_2$, $R(e_3) = 0$,

(RB6) $R(e_1) = e_2$, $R(e_2) = R(e_3) = 0$,

(RB7) $R(e_1) = e_2$, $R(e_2) = 0$, $R(e_3) = -e_3$.

{\sc Proof}. 
a) Non-splitting case appears only when the graph $T'$ is non-empty and improperly 2-colored.
Up to $\phi$, we may assume that two vertices are colored in white.

b) Cases (RB1)--(RB3) correspond to improperly 2-colorings of the graph $T'$
with $V(T') = \{1,2,3\}$ and $E(T') = \{(1,2),(2,3)\}$.
Cases (RB4), (RB5) correspond to improperly 2-colorings of the graph $T'$
with $E(T') = \{(1,2),(1,3)\}$. Finally, cases (RB6), (RB7) correspond to 
improperly 2-colorings of the graph $T'$ with $E(T') = \{(1,2)\}$.

{\bf Statement 6}.
Up to $\phi$, we have the following splitting but not inner-splitting 
pairwise nonisomorphic RB-operators:

a) for $n=2$: $R(e_1) = -e_1$, $R(e_2) = 0$;

b) for $n=3$:

(RB1$'$) $R(e_1) = e_2$, $R(e_2) = 0$, $R(e_3) = -e_3$,

(RB2$'$) $R(e_1) = -e_1$, $R(e_2) = R(e_3) = 0$.

\section{RB-induced algebra structures on a sum of fields}

Let $C$ be an associative algebra and $R$ be an RB-operator on $C$ of weight~$\lambda$.
Then the space~$C$ under the product 
\begin{equation}\label{postProduct}
x\circ_R y  = R(x)y + xR(y) + \lambda xy
\end{equation}
is an associative algebra~\cite{GuoMonograph,Embedding}. 
Let us denote the obtained algebra as~$C^R$.
It is easy to see that $C^{\phi(R)}\cong C^R$.

Let us denote by $\mathrm{Ab}_n$ the $n$-dimensional algebra
with zero (trivial) product.

{\bf Theorem 3}.
Given an algebra $A = Fe_1\oplus \ldots \oplus Fe_n$ 
and an RB-operator $R$ of weight~$\lambda$ on $A$, we have 
$A^R \cong \begin{cases} \mathrm{Ab}_n, & \lambda = 0, \\
A, & \lambda\neq0. \end{cases}$

{\sc Proof}.
If $\lambda = 0$, then $R = 0$~\cite{Unital} and $x\circ_R y = 0$.
For $\lambda\neq0$, we may assume that $\lambda = 1$, 
since rescalling of the product does not exchange the algebraic structure.

Let us prove the statement by induction on~$n$.
For $n=1$, we have either $R = 0$ or $R = -\id$. 
Due to~\eqref{postProduct} we get either $x\circ y = xy$ or $x\circ y = -xy$,
in both cases $A^R\cong A$.

Suppose that we have proved Theorem~3 for all numbers less $n$.
Let a~graph $T' = T'(R)$ with $n$ vertices be disconnected, 
denote by $T_1,\ldots,T_k$ the connected components of $T'$.
As earlier, we define $A = A_1\oplus \ldots A_k$ for $A_s = \Span\{e_j\mid j\in V(T_s)\}$
and define $R_s$ as the induced RB-operator $R|_{A_s}$.
By the induction hypothesis, $A_s^R\cong A_s$ for every $s$ and so
$A = A_1\oplus \ldots \oplus A_k\cong A_1^R\oplus \ldots \oplus A_k^R = A^R$.

Now consider the case when $T'$ is connected. 
We may assume that $e_1$ corresponds to the vertex~1, the only source in $G$.
Note the space $I_1 = \Span\{e_j\mid j\geq2\}$ is an ideal in $A^R$
which is isomorphic to $Fe_2\oplus \ldots \oplus Fe_n$ by the induction hypothesis.
Up to $\phi$, we may assume that~the vertex~1 in $T'$ is colored in white
and $2,\ldots,t$ is a list of all neighbours of~1 in $T'$.
Let us consider the one-dimensional space $I_2$ in $A^R$ generated by the vector 
$a = e_1 - c(2)e_2 - \ldots - c(t)e_t$, where 
$$
c(i) = \begin{cases} 
1, & i\,\mbox{ is colored in white}, \\
-1, & i\,\mbox{ is colored in black}.
\end{cases}$$ 
In terms of the matrix entries, $c(i) = 1+2r_{ii}$.
We may assume that $c(2) = c(3) = \ldots = c(s) = 1$ and 
$c(s+1) = \ldots = c(t) = -1$ for some $s\in\{2,\ldots,t\}$. 

By~\eqref{postProduct} we compute the product of $a$ with $e_k$ for $k>t$:
\begin{multline*}
a\circ e_k 
 = (e_1+e_2+\ldots+e_s-e_{s+1}-\ldots-e_t)\circ e_k \\
 = R(e_1+e_2+\ldots+e_s-e_{s+1}-\ldots-e_t)e_k. 
\end{multline*}
Since $k$ is connected with only one vertex from $2,\ldots,t$ (due to (SF2)),
say $j$, we have 
$$
a\circ e_k 
 = R(e_1 - c(j)e_j)e_k
 = e_k - c(j)(1+2r_{jj})e_k 
 = (1-(c(j))^2)e_k
 = 0.
$$
Analogously we can check that $a\circ e_k = 0$ for all $k>1$.
Thus, $I_2$ is an ideal in $A^R$. 

Now, we calculate 
\begin{multline*}
a\circ a 
    = e_1\circ (e_1+e_2+\ldots+e_s-e_{s+1}-\ldots-e_t) \\
    = R(e_1)(e_1+e_2+\ldots+e_s-e_{s+1}-\ldots-e_t) + e_1 \\
    = (e_2+\ldots+e_s+e_{s+1}+\ldots+e_t)(e_1+e_2+\ldots+e_s-e_{s+1}-\ldots-e_t) + e_1 \\
    = e_1 + e_2 + \ldots + e_s - e_{s+1} - \ldots -e_t = a
\end{multline*}
and so $I_2$ is isomorphic to~$F$.

Summarising, we have
$A^R = I_1\oplus I_2 \cong (Fe_2\oplus \ldots \oplus Fe_n)\oplus F \cong A$.
Theorem is proved.

\section*{Acknowledgements}

The main part of the paper was done while working in Sobolev Institute of Mathematics in 2017.
The research is supported by RSF (project N 14-21-00065).

\noindent Vsevolod Gubarev \\
University of Vienna \\
Oskar-Morgenstern-Platz 1, 1090 Vienna, Austria \\
Sobolev Institute of mathematics \\
Acad. Koptyug ave. 4, 630090 Novosibirsk, Russia \\
e-mail: vsevolod.gubarev@univie.ac.at

\begin{thebibliography}{67}
\bibitem{AnBai}
H. An, C. Bai.
From Rota-Baxter Algebras to Pre-Lie Algebras.
J. Phys. A (1) (2008), 015201, 19~p.

\bibitem{Atkinson}
F.V. Atkinson.
Some aspects of Baxter’s functional equation.
J. Math. Anal. Appl. {\bf 7} (1963) 1-30.

\bibitem{Baxter}
G. Baxter. 
An analytic problem whose solution follows from a simple algebraic identity. 
Pacific J. Math. {\bf 10} (1960) 731--742.

\bibitem{BelaDrin82}
A.A. Belavin, V.G.  Drinfel'd. 
Solutions of the classical Yang---Baxter equation for simple Lie algebras.
Funct. Anal. Appl. (3) {\bf 16} (1982) 159--180.

\bibitem{BGP}
P. Benito, V. Gubarev, A. Pozhidaev.
Rota---Baxter operators on quadratic algebras.
Mediterr. J. Math. {\bf 15} (2018), 23~p. (N189). 

\bibitem{Braga}
S.L. de Bragan\c{c}a. 
Finite Dimensional Baxter Algebras. Stud. Appl. Math. (1) {\bf 54} (1975) 75--89.

\bibitem{BurdeDekimpeVercammen12} 
D. Burde, K. Dekimpe and K. Vercammen. 
Affine actions on Lie groups and post-Lie algebra structures.
Linear Algebra Appl. (5) {\bf 437} (2012) 1250--1263.

\bibitem{BurdeDekimpe13} 
D. Burde, K. Dekimpe. 
Post-Lie algebra structures and generalized derivations of semisimple Lie algebras. 
Mosc. Math. J. (1) {\bf 13} (2013) 1--18.

\bibitem{BurdeGub18} 
D. Burde, V. Gubarev. 
Rota---Baxter operators and post-Lie algebra structures on semisimple Lie algebras.
Commun. Algebra (accepted), arXiv:1805.05104 [RA], 18~p.

\bibitem{Cartier}
P. Cartier.
On the structure of free Baxter algebras. Adv. Math. {\bf 9} (1972) 253--265.

\bibitem{Castelo}
R. Castelo and A. Siebes. 
A characterization of moral transitive acyclic directed graph Markov models as labeled trees. 
J. Stat. Plan. Inf. {\bf 115} (2003) 235--259.

\bibitem{Unital}
V. Gubarev. 
Rota---Baxter operators on unital algebras. arXiv.1805.00723v2, 37~p.

\bibitem{Embedding}
V. Gubarev, P. Kolesnikov.
Embedding of dendriform algebras into Rota---Baxter algebras.
Cent. Eur. J. Math. (2) {\bf 11} (2013) 226--245.

\bibitem{GuoMonograph}
L. Guo.
An Introduction to Rota---Baxter Algebra. Surveys of Modern Mathematics, vol. 4,
Int. Press, Somerville (MA, USA); Higher education press, Beijing, 2012.

\bibitem{Harary}
F. Harary, J. Kabell, F. McMorris. 
Subtree acyclic digraphs. Ars Combin. {\bf 34} (1992) 93--95.

\bibitem{Lauritzen} 
S. Lauritzen. 
Graphical Models. Oxford: Oxford University Press, 1996.

\bibitem{Miller}
J.B. Miller. 
Baxter operators and endomorphisms on Banach algebras.
J. Math. Anal. Appl. {\bf 25} (1969) 503--520.

\bibitem{OEIS}
OEIS Foundation Inc. The on-line encyclopedia of integer sequences, http://oeis.org.

\bibitem{Rota}
G.-C. Rota. 
Baxter algebras and combinatorial identities. I. 
Bull. Amer. Math. Soc. {\bf 75} (1969) 325--329.

\bibitem{Semenov83}
M.A. Semenov-Tyan-Shanskii. 
What is a classical $r$-matrix? 
Funct. Anal. Appl. {\bf 17} (1983) 259--272.

\end{thebibliography}
\end{document}